\documentclass[12pt]{amsart}
\usepackage{amssymb,latexsym}

\numberwithin{equation}{section}
\newtheorem{thm}{Theorem}[section]

\newtheorem{prop}[thm]{Proposition}
\newtheorem{cor}[thm]{Corollary}
\newtheorem{rem}[thm]{Remark}
\newtheorem{lem}[thm]{Lemma}
\newtheorem{example}[thm]{Example}
\newtheorem{definition}[thm]{Definition}

\begin{document}

\title{ALGEBRAS ASSOCIATED TO  ACYCLIC DIRECTED GRAPHS}

\author{Vladimir Retakh}
\author{Robert Lee Wilson}
\address{Department of Mathematics, Rutgers University, Piscataway, NJ 08854-8019, USA}
\email{vretakh@math.rutgers.edu}
\email{rwilson@math.rutgers.edu}
\thanks{Both authors were supported in part by
the NSA grant H98230-06-1-0028}

\keywords{generalized layered graphs, Hilbert series, factorizations of noncommutative polynomials}
\subjclass{05E05; 15A15; 16W30}

\begin{abstract}
We construct and study a class of algebras associated to generalized layered graphs,
i.e. directed graphs with a ranking function on their vertices and edges.
Each finite acyclic directed  graph admits countably many  structures
of a  generalized layered graph. We construct linear bases in such algebras and compute their Hilbert series.
Our interest to generalized layered graphs and algebras associated
to those graphs is motivated by their relations to factorizations
of polynomials over noncommutative rings.
\end{abstract}

\maketitle

\section{Introduction}

By a generalized layered graph we mean a pair $\Gamma = (G,|.|)$ 
where $G = (V, E)$ is a directed graph and  
$|.|:V \rightarrow {\mathbf Z}_{\ge 0}$ satisfies $|v| > |w|$ 
whenever $v,w \in V$ and there is an edge $e \in E$ from $v$ to $w$. We call $|.|$ the rank function of $\Gamma$. 
We write $l(e) = |v| - |w|$ and call this the length of the edge $e$.  We will
see that if $G$ is any acyclic directed graph then there are 
countably many rank functions $|.|$ such that $(G,|.|)$ is a generalized layered graph.

In this paper we construct and study a class of algebras $A(\Gamma)$ associated to generalized
layered graphs $\Gamma$. Generators of our algebras are
elements $a_1(e), a_2(e),\dots , a_{l(e)}(e)$ associated to edges $e$ of $\Gamma$.
The relations are defined as follows. Let sequences of edges
$e_1, e_2, \dots , e_p$ and $f_1,f_2,\dots , f_q$ define paths
with the same end and the same origin. Then they define a
relation given by the identity
$$U_{e_1}(\tau)U_{e_2}(\tau)\dots U_{e_p}(\tau)=U_{f_1}(\tau)U_{f_2}(\tau)\dots U_{f_q}(\tau),$$
where $\tau $ is a formal central variable and
$$U_e(\tau)=\tau ^{l(e)}- a_1(e)\tau ^{l(e)-1}+ a_2(e)\tau ^{l(e)-2}-\dots \pm a_{l(e)}(e)$$
for  any edge $e$ of $\Gamma$.  We will show that if $\Gamma = (G,|.|)$ the 
structure of $A(\Gamma)$ depends on the rank function $|.|$ as well as directed graph $G$.

Our interest to generalized layered graphs and algebras associated
to those graphs is motivated by their relations to factorizations
of polynomials over noncommutative rings. Let $R$ be a unital
algebra, $P(\tau)$ a monic polynomial over $R$, and  $\mathcal P$ be
a set of monic right divisors of $P(\tau)$, i.e.  $\mathcal P$
consists of monic polynomials $Q(\tau)\in R[\tau]$ such that
$P(\tau)=U(\tau)Q(\tau)$ for some $U(\tau)\in R[\tau]$. In papers
\cite{GRW, GRSW, GGRW, RSW3} we studied subalgebras $R_{\mathcal P}$ of $R$ generated by the
coefficients of polynomials $U(\tau)$ for certain sets $\mathcal P$.
These studies led us to our notion of universal algebras of type
$R_{\mathcal P}$ defined via generalized layered graphs.

To any set $\mathcal P$ we associate a generalized  layered graph $\Gamma
(\mathcal P)$ constructed in the following way. The vertices of rank $m$
in $\Gamma $ are polynomials  $Q(\tau)\in \mathcal P$ of degree $m$. An edge $e$
of  length $k$ goes from a vertex $Q_1(\tau)$ to a vertex $Q_2(\tau)$ if and only if
 $Q_1(\tau)=U(\tau)Q_2(\tau)$ where $U(\tau)\in R[\tau]$ is a polynomial
of degree $k$.  A set of polynomials $Q_j(\tau)\in \mathcal P$, $j=1,2,\dots ,p$ such that
$Q_j(\tau)=U_j(\tau)Q_{j+1}(\tau)$ for $j=1,2,\dots , p-1$ defines
a path in  $\Gamma (\mathcal P)$  from vertex $Q_1(\tau)$ to a vertex $Q_p(\tau)$. If a set
of polynomials $S_k(\tau)$, $k=1,2,\dots , q$ defines another path (say $S_k(\tau)=W_k(\tau)S_{k+1}(\tau)$)
in this graph
with the same origin (i.e.  $Q_1(\tau)=S_1(\tau)$) and the same end (i.e. $Q_p(\tau)=S_q(\tau)$)
then $U_1(\tau)U_2(\tau)\dots U_{p-1}(\tau)= W_1(\tau)W_2(\tau)\dots W_{q-1}(\tau)$.

If $\mathcal P$ contains polynomials $1$ and $P(\tau)$ then the
graph $\Gamma (\mathcal P)$ contains exactly one minimal vertex
$1\in R$ of rank $0$ and exactly one maximal vertex $P(\tau)$ of
rank $\deg P(\tau)$.

 Note that there is a canonical homomorphism $\phi: A(\Gamma (\mathcal P))\rightarrow R$ defined in the
following way. For any edge $e$ in $\Gamma (\mathcal P)$ of  length
$k$ there are corresponding generators $a_1(e), a_2(e),\dots ,
a_k(e)$ in $A(\Gamma (\mathcal P))$.  Since the edge goes from a vertex associated to some polynomial 
$Q(\tau)$ to a vertex associated to some polynomial $W(\tau)$  there is a polynomial
$$S_e(\tau)=\tau^k-b_1\tau^{k-1}+b_2\tau^{k-2}-\dots
+(-1)^kb_k \in R[\tau]$$ 
such that $Q(\tau) = S_e(\tau)W(\tau).$ 
Then there is a unique $\phi$ satisfying $\phi (a_i(e))=b_i$, $i=1,2,\dots ,k.$
The image of $\phi $ is a subalgebra  $R_{\mathcal P}$ of $R$ which 
records information about factorizations of the initial
polynomial $P(\tau)$.

This paper continues our investigations started in \cite{GRSW, GGRW, RSW1, RSW2, RSW3}  where we defined
and studied properties of
the algebras $A(\Gamma )$ for layered graphs $\Gamma $, i.e. graphs with edges of length one. These algebras
correspond to factorizations of polynomials into products of linear factors. The algebras $A(\Gamma )$
for layered graphs have a deep interesting structure.  In this paper we show that many results about the  algebras
$A(\Gamma )$ for a layered graph $\Gamma $
admit natural generalization to the much wider class of generalized layered graphs.

For reasons of clarity we have required that the polynomials in the foregoing discussion be monic.  
However, in the body of the paper, in order to be able to write inverses easily in terms of 
geometric series, we will work instead with polynomials with constant term $1$.

The paper is organized as follows. In Section 1 we prove that any acyclic directed graph may be given the structure
of a generalized layered graph in countably many ways and define the algebras $A(\Gamma)$ for generalized layered graphs.
In Section 2 we construct a spanning set of monomials in the algebra $A(\Gamma)$. In Section 3 we prove that this
spanning set is, in fact, a linear basis in $A(\Gamma)$. Our proof is based on the main result from \cite{GRSW}
and an analysis of the behavior of algebras $A(\Gamma)$ under a natural operation on $\Gamma$ (adding a vertex).
We also study the behavior of these algebras under other operations (adding edges, etc.) and show that our
construction gives a functor from the category of generalized layered graphs to the category of
associative algebras. 
In Section 4 we compute Hilbert series for the algebras $A(\Gamma)$. This generalize the
main result for layered graphs from \cite{RSW2}. In Section 5 we study the behavior of the Hilbert series under operations
on graphs. In particular, we use these results to construct a natural family of noncommutative complete intersections 
(in the sense of \cite{Anick}, \cite{Golod}.) 


\section{Basic definitions}

Let $F$ be a field and for any set $S$ let $T(S)$
denote the free associative algebra on $S$ over $F$.

Let $G = (V,E)$ be a directed graph.  Thus $V$ is a set
whose elements are called the {\it vertices} of $G$,
$E$ is a set whose elements are called the {\it edges}
of $G$, and there are two functions $h,t: E \rightarrow V$.
For $e \in E, t(e)$ is called the {\it tail} of $e$
and $h(e)$ is called the {\it head} of $e$.
We say that a pair $\Gamma = (G,|.|)$ is a {\it generalized layered graph} 
if $G =(V,E)$ is a directed graph and $|.|:V \rightarrow {\bold Z}_{\ge 0}$ 
satisfies $|v| > |w|$ whenever there is an edge $e$ from $v$ to $w$.  
If $\Gamma$ is a generalized layered graph and $V_i$ is the set of $v\in V$
such that $|v| = i$, 
then $V$ is the disjoint union $\cup_{i=0}^{\infty} V_i.$
If $e \in E$ we define $|e| = |t(e)|$ and $l(e)$,
the {\it length} of $e$, to be $|t(e)| - |h(e)|$.
Recall that if $l(e) = 1$ for all $e \in E$, then $\Gamma$ is a {\it layered graph}.

While the requirement that a directed graph is layered is quite restrictive,
the requirement that a graph be generalized layered is not.

\begin{prop}
Let $G= (V,E)$ be a finite acyclic directed graph.  Then there are countably many ranking functions
$|.|$ such that $(G,|.|)$ is a generalized layered graph.
\end{prop}

\begin{proof}  The result clearly holds if $|V| = 1$.  We will proceed by induction on $|V|$
Let $v' \in V$ satisfy $\{e \in E|h(e) = v'\} = \emptyset.$  Set $V' = V \setminus \{v'\}$ and
$E' = \{e \in E|t(e) \ne v'\}$.  Then $(V',E')$ is a finite acyclic 
directed graph and so, by the induction assumption, has countably many ranking functions.  
Let $|.|'$ be one such
and let 
$|.|:V \rightarrow {\bold Z}_{\ge 0}$ extend $|.|'$.  Then $|.|$ 
is a ranking function on $G$ if and only if 
$|v'|> |w|$ for all $w \in V'$ such that there is an edge from $v'$ to $w$.  
Thus $|.|'$ has countably 
many extensions to ranking functions of $G$.  
\end{proof}

\begin{rem} The above argument shows that if $|v|_{can}$ is defined to be $0$ 
when $\{e\in E|t(e) = v\} = \emptyset$ and otherwise defined to be
the largest $s$ such that there exist edges $e_1,...,e_s \in E$ with $t(e_1) = v$ 
and $h(e_i) = t(e_{i+1})$ for $1 \le i < s.$,
then $(G,|.|_{can})$ is a generalized layered graph and if 
$(G,|.|)$ is any generalized layered graph, then $|v| \ge |v|_{can}$ for all $v \in V.$
\end{rem}

Let $G =(V,E)$ be a acyclic directed graph. Let ${\mathcal G}(G)$ denote 
the (countable) set of all generalized layered graphs
$(G,|.|)$. For $(G,|.|_1)$, $(G,|.|_2) \in {\mathcal G}(G)$ we write
$(G,|.|_1) \ge (G,|.|_2)$ if $|e|_1 \ge |e|_2$ for all $e \in E$.  
With this definition, ${\mathcal G}(G)$ becomes a partially ordered set. 

We will now define the algebras which are the primary objects of study in this paper.
If $G = (V,E)$ and $\Gamma = (G,|.|)$ 
is a generalized layered graph, let $E^{\sharp}$ denote the subset
$$\cup_{e \in E} \{e\} \times [1,l(e)]$$
of $E \times {\mathbf Z}$.  Denote the ordered pair
$(e,i) \in E^{\sharp}$ by $a_i(e)$.  For $a_i(e) \in E^{\sharp}$
define $\deg (a_i(e)) = i.$ Let $t$ be a central variable and define 
$P_e(t)$ to be the polynomial
$$1 + \sum_{j=1}^{l(e)} (-1)^j a_j(e)t^j$$
in $T(E^{\sharp})[t]$.   Set $a_0(e) = 1$.
Then $deg$ gives $T(E^{\sharp})$ the structure of a graded algebra
$$T(E^{\sharp}) = \oplus T(E^{\sharp})_{[i]}.$$
Note that if we set $\deg(t) = -1$ then
$P_e(t)$ is homogeneous of degree $0$ in the graded algebra $T(E^{\sharp})[t]$.

Furthermore, we may also give $T(E^{\sharp})$ the structure of a filtered algebra.
Set $$f(a,j)=ja-j(j-1)/2$$ and for $(e,i)\in E^{\sharp}$ set 
$|a_i(e)|=\sum _{j=1}^i(|e|-j+1)=f(|e|,i)$. Then setting
$$T(E^{\sharp})_i = span\{a_{j_1}(e_1)...a_{j_r}(e_r)|r \ge 0, |a_{j_1}(e_1)|  + ... + |a_{j_r}(e_r)| \le i\}$$
gives $T(E^{\sharp})$ the structure of a filtered algebra.

As usual, we say that a sequence $\pi = (e_1,...,e_r)$ of
edges of $\Gamma$ is a {\it (directed) path} from $v$ to $w$ if
$t(e_1) = v, h(e_r) = w$ and
$h(e_j) = t(e_{j+1})$ for all $j, 1 \le j < r$.
If $\pi$ is such a path, we write $v > w$, define
$t(\pi) = t(e_1)$ and call this the tail of $\pi$, and
define $h(\pi) = h(e_r)$ and call this the head of
$\pi$. We define $|\pi| = |t(\pi)|$ and $l(\pi) = |t(e_1)| - |h(e_r)| =
\sum_{i=1}^r l(e_i)$.  For any path $\pi = (e_1,...,e_r)$ in $\Gamma$ define
$$P_{\pi}(t) = P_{e_1}(t)P_{e_2}(t)...P_{e_r}(t) \in T(E^{\sharp})[t]$$
and write
$$P_{\pi}(t) = 1 + \sum_{j=1}^{l(\pi)} (-1)^je(\pi,j)t^j.$$
Let $R$ denote the ideal in $T(E^{\sharp})$ generated by
$$\{e(\pi_1,j) - e(\pi_2,j)\ |\ t(\pi_1) = t(\pi_2), h(\pi_1) = h(\pi_2),\ 1\leq j\leq l(\pi_1)\}.$$

\begin{definition} For a generalized layered graph $\Gamma = (G, |.|)$ where $G = (V,E)$ we define
$$ A(\Gamma) = T(E^{\sharp})/R.$$
\end{definition}

Thus the images of $P_{\pi_1}(t)$ and $P_{\pi_2}(t)$ in $A(\Gamma)[t]$ are equal whenever
$t(\pi_1) = t(\pi_2), h(\pi_1) = h(\pi_2)$.

Since each $P_{\pi}(t)$ is homogeneous of degree $0$, the ideal $R$ is
generated by homogeneous elements and so is a graded ideal.  Thus
$$A(\Gamma) = \oplus A(\Gamma)_{[i]}$$
has the structure of a graded algebra where
$$A(\Gamma)_{[i]} = (T(E^{\sharp})_{[i]} + R)/R.$$
$A(\Gamma)$ also has the structure of a filtered algebra where
$$A(\Gamma)_i = (T(E^{\sharp})_i + R)/R.$$

When the meaning is clear from context we will denote elements of $A(\Gamma)$ or 
$A(\Gamma)[t]$ by their representatives in $T(E^{\sharp})$ or $T(E^{\sharp})[t].$

Recall that $f(a,j) = ja - j(j-1)/2.$

\begin{lem} $e(\pi,j) \in A(\Gamma)_{f(|\pi|,j)}$
for $1 \le j \le l(\pi)$.
\end{lem}

\begin{proof}  Write $\pi = (e_1,...,e_r)$.
Let $U$ denote the set of all sequences of
integers $(i_1,...,i_r)$ with $i_1 + ... + i_r = j$
and $0 \le i_k \le l(e_k)$ for $1 \le k \le r$.
Then
$$e(\pi,j) = \sum_{(i_1,...,i_r) \in U} a_{i_1}(e_1)...a_{i_r}(e_r).$$
For $(i_1,...,i_r) \in U$, write 
$$S(i_1,...,i_r) = \sum_{k=1}^r (|e_k| + (|e_k| - 1 ) + ... + (|e_k| - i_k + 1)).$$
Then $$a_{i_1}(e_1)...a_{i_r}(e_r) \in A(\Gamma)_{S(i_1,...,i_r)}.$$
Now if $(i_1,...,i_r), (i_1,...,i_{k-1},i_k + 1, i_{k+1} - 1, i_{k+2},...,i_r) \in U$
we have
$S(i_1,...,i_r) \le S(i_1,...,i_{k-1},i_k + 1, i_{k+1} - 1, i_{k+2},...,i_r)$.
Then writing
$j = l(e_1) + l(e_2) + ... + l(e_s) + j'$ where $0 \le s < r$ and $0 \le j' \leq l(e_{s+1})$
we see that
for any $(i_1,...,i_r) \in U$ we have
$$S(i_1,...,i_r) \le S(l(e_1),...,l(e_s),j',0,...,0).$$
Since $S(l(e_1),...,l(e_s),j',0,...,0) = j|\pi|-(j(j-1)/2) = f(|\pi|,j)$ the lemma is proved.
\end{proof}

We may also define a sequence of related algebras, $A(k,\Gamma), k \ge 1$,
by requiring that the images of $P_{\pi_1}(t)$
and $P_{\pi_2}(t)$ in $A(\Gamma)[t]/(t^k)$ are equal whenever
$t(\pi_1) = t(\pi_2), h(\pi_1) = h(\pi_2)$.
Thus we define $R_k$ to be the ideal in $T(E^{\sharp})$ generated by
$$\{e(\pi_1,j) - e(\pi_2,j)|1 \le j < k, t(\pi_1) = t(\pi_2), h(\pi_1) = h(\pi_2)\}$$
and define $$A(k,\Gamma) = T(E^{\sharp})/R_k.$$
Then $$A(1,\Gamma) \rightarrow A(2,\Gamma) \rightarrow ...$$
and if $V = \coprod_{j=0}^n V_j$ we have $A(\Gamma) = A(m,\Gamma)$ whenever $m > n$.
Note that the image of $P_{\pi}(t)$ in $T(E^{\sharp})[t]/(t^j)$ is
invertible for any path $\pi$ and for any $j > 0$ and that its
inverse is given by the image of the geometric series:
$$P_{\pi}(t)^{-1} = \sum_{m=0}^{j-1} (1 - P_{\pi}(t))^m.$$
\medskip
\begin{rem} One might require $P_e(t)$ to be a
monic polynomial (instead of requiring that the constant term be
$1$). This gives an equivalent definition.  We have chosen to
require that the constant term be $1$ in order to have the above
easy expressions for inverses.
\end{rem}
\begin{example} Let $\Gamma = (G,|.|)$ and assume $G = (V,E)$ is acyclic as a (nondirected) graph, i.e., a tree. 
Let $E = \{e_1,...,e_n\}$ and recall that
$l(e)$ denotes the length of an edge $e$. Write $K = \sum_{i=1}^n
l(e_i)$ and assume that the edges have been ordered so that
$l(e_i) \le l(e_{i+1})$ for $1 \le i < n$ so that $(l_1,...,l_n)$
is a partition of $K$. Let $(m_1,...,m_r)$ be the partition
conjugate to $(l_1,...,l_n).$ Then $R = (0)$ and so $A(\Gamma) = T(E^{\sharp})$ is the free algebra on 
$K$ generators, $m_j$ of which have degree $j$.
\end{example}
\medskip
Algebras associated with generalized layered graphs occur naturally in the study of certain modules for the symmetric group.  
The following example is taken from Duffy, \cite{D1,D2}.
\begin{example}  Let $[n]$ denote the
set $\{1,...,n\}$. For $ \sigma \in Sym_n$ define $[n,\sigma]$ to
be the set of $\sigma$-orbits in $[n]$. Let ${\mathcal
P}([n,\sigma])$ denote the power set of $[n,\sigma]$ and
$G([n,\sigma])$ denote the Hasse graph of the partially
ordered set ${\mathcal P}([n,\sigma]).$ For $v = \{v_1,...,v_k\}
\in {\mathcal P}([n,\sigma])$ let $|v|$ denote the cardinality of
$\cup_{i=1}^k v_i$.  Then $\Gamma([n,\sigma]) = (G([n,\sigma]),|.|)$ is a generalized
layered graph. In particular, $\Gamma([3,(12)])$ has four
vertices: $a = \{\{1,2\},\{3\}\}, b = \{\{1,2\}\}, c = \{\{3\}\},
* = \emptyset$ and edges $e_1$ from $a$ to $b$, $e_2$ from $a$ to
$c$, $e_3$ from $b$ to $*$, and $e_4$ from $c$ to $*$. Here $|a|
= 3, |b| = 2, |c| = 1, |*| = 0, l(e_1) = l(e_4) = 1, l(e_2) =
l(e_3) = 2.$
\end{example}

We define now the category of generalized layered graphs.
Let $\Gamma = (G,|.|)$ and $\Gamma' = (G',|.|')$ be generalized layered graphs where 
$G = (V,E)$ and $G' = (V',E')$ are directed graphs.  A morphism $\phi:\Gamma \rightarrow \Gamma'$ is a pair
$\phi = (\phi_V,\phi_E)$ where $$\phi_V:V \rightarrow V',$$ 
$$ \phi_E:E \rightarrow E',$$ 
such that for all $e\in E$
$$t(\phi_E(e)) = \phi_V(t(e)),$$
$$h(\phi_E(e)) = \phi_V(h(e)), $$
and
$$l(\phi_E(e)) \le l(e).$$

We denote by ${\mathcal {GLG}}$ the category with generalized layered graphs as objects and with morphisms as defined above.

Let $\phi \in \hom(\Gamma,\Gamma')$ where $\Gamma = (G,|.|), \Gamma' = (G',|.|')$.  For a ground field $F$
define
$$\tilde\phi_E:E^{\sharp} \rightarrow FE^{'\sharp}$$
by
$$\tilde\phi_E(a_i(e)) = a_i(\phi_E(e))$$
if $0 \le i \le l(\phi_E(e)),$ and
$$\tilde\phi_E(a_i(e)) = 0$$
if $i > l(\phi_E(e)).$

Then $\tilde\phi_E$ induces a homomorphism, again denoted $\tilde\phi_E$ from
$T(E^{\sharp}) \rightarrow T(E^{'\sharp})$
and hence also induces a homomorphism, still denoted $\tilde\phi_E$, from
$T(E^{\sharp})[t] \rightarrow T(E^{'\sharp})[t].$
Clearly
$\tilde\phi_E(P_e)(t) = P_{\phi_E(e)}(t)$
and so, if $\pi $ is a path,
$\tilde\phi_E(P_{\pi}(t)) = P_{\phi_E(\pi)}(t)$.  Therefore
$\tilde\phi(R_{\Gamma}) \subseteq R_{\Gamma'}$ and hence 
$\tilde\phi$ induces a homomorphism of associative algebras
$$A(\phi):A(\Gamma) \rightarrow A(\Gamma').$$

\begin{prop} The mappings
$$\Gamma \mapsto A(\Gamma)$$
and
$$\phi \mapsto A(\phi)$$
define a functor from $\mathcal{GLG}$ to the category of associative algebras.
\end{prop}
Addition of a single edge gives an interesting example of a morphism of 
generalized layered graphs.  Thus let $G = (V,E)$ be a directed graph and $v,w \in V, |v| > |w|.$  
Define $G^e  = (V,E^e)$ by   $E^e = E \cup \{e\}, t(e) = v, h(e) = w.$  Let $|.|$ be
a ranking function for $G$ and $\Gamma = (G,|.|), \Gamma^e = (G',|.|).$  
Let $i_V$ denote the identity map on $V$ and $i^e_E$ denote the 
injection of $E$ into $E^e$.  Then $i^e = (i_V,i^e_E)$ is a morphism of $\Gamma$ to $\Gamma^e$. 
We will consider properties of the corresponding homomorphism $A(i^e)$ in Section 3.

\section{Spanning set for $A(\Gamma)$}

Let $\Gamma = (G,|.|)$ where $G = (V,E)$ be a generalized layered graph and
assume that $V_0$ is the set of minimal vertices of $\Gamma$.
Write $V_+ = V \setminus V_0$.

For each vertex $v \in V_+$ fix (arbitrarily) an edge
$e_v \in E$ with $t(e_v) = v$.  Then there is a path $\pi_v$ defined by
$\pi_v = (e_1,...,e_r)$ with $e_1 = e_v, e_{i+1} = e_{h(e_i)}$ for $1 \le i <r$,
$h(e_r) \in V_0$.  Set $P_v(t) = P_{\pi_v}(t)$. Let $V^{\sharp}$ denote the subset
$$\cup_{v \in V_+} \{v\} \times [1,|v|]$$
of $V_+ \times {\mathbf Z}$  and, for $(v,j) \in V^{\sharp} $ set $e(v,j) = e(\pi_v,j)$.
Define a partial order on $V \times {\mathbf Z}$ by
$(v,k) \gtrdot (w,l)$ if $v > w$ and $k = |v| - |w|$.
\medskip
\begin{lem}    Let $\pi$ be a path with $t(\pi) = v$ and $l(\pi) = m$.
$$e(\pi,m) \equiv e(v,m) \pmod {A(\Gamma)_{f(|v|,m)-1}}.$$
\end{lem}

\begin{proof}  Set $w = h(\pi)$.  Then $P_{\pi}(t) = P_v(t)P_w(t)^{-1}$ 
in $A(\Gamma )[[t]]$ and so (since the constant term of $P_w(t)$ is $1$)
$$e(\pi,m) = \sum (-1)^r e(v,i_0)e(w,i_1)...e(w,i_r)$$
where the sum is taken over all sequences of integers
$(i_0,...,i_r)$ with $r \ge 0, i_0 \ge 0, i_1,...,i_r \ge 1$ and $i_0 + ... + i_r = m$.
Let $$M(i_0,...,i_r) =  |v| + (|v|-1) + ... + (|v| -i_0 + 1) + $$
$$\sum_{j=1}^r (|w| + (|w| - 1|) + ...  + (|w| - i_j + 1)).$$
\medskip
Clearly $M(i_0,...,i_r) \le M(i_0) + (m-i_0)|w|$ and as $|v|-m=|w|$
$M(i_0) + (m-i_0)|w| \le M(m)$ with equality if and only if $i_0 = m$
(or, equivalently, $r = 0$). Since $e(v,i_0)e(w,i_1)...e(w,i_r) \in A(\Gamma)_{M(i_0,...,i_r)}$ 
and $M(m) = f(|v|,m)$, the lemma is proved.
\end{proof}
\medskip
\begin{prop} If $(v,k) \gtrdot (w,l)$ then
$$e(v,k)e(w,l) \equiv e(v,k+l) \pmod {A(\Gamma)_{f(|v|,k+l)-1}}.$$
\end{prop}

\begin{proof}  Let $\pi$ be a path from $v$ to $w$.
Then $P_v(t) = P_{\pi}(t)P_w(t).$
Since $l(\pi) = |v| - |w| = k$, $P_{\pi}(t)$ is a polynomial
of degree $k$ and so $e(v,k+l)$, the coefficient of $(-t)^{k+l}$ in $P_v(t)$, is equal to
$$\sum_{j=0}^k e(\pi,j)e(w,k+l-j).$$
By Lemma 1.4,
$$e(\pi,j)e(w,k+l-j) \in A(\Gamma)_{f(|\pi|,j) + f(|w|,k+l-j)}.$$
Now, as $0\le j\le k$, $$f(|\pi|,j) + f(|w|,k+l-j) \le f(|v|,k+l)$$ with equality if
and only if $j = k+l$. Thus $$e(v,k+l) \equiv e(\pi,k)e(w,l)
\pmod {A(\Gamma)_{f(|v|,k+l) - 1}}$$ and so the proposition
follows from Lemma 2.1.
\end{proof}

Define $$\mathbf B_1(\Gamma) = \{((v_1,k_k),...,(v_r,k_r))\ |\ r \ge 0,
v_1,...,v_r \in V_+, 1 \le k_i \le |v_i|\}$$ and 
$$\mathbf B(\Gamma) =
\{((v_1,k_1)...(v_r,k_r) \in \mathbf B_1(\Gamma)\ |\ (v_{j},i_{j})
\not\gtrdot (v_{j+1},i_{j+1}), 1 \le j < r\}.$$ 
Define
$$\epsilon:\mathbf B_1(\Gamma) \rightarrow A(\Gamma)$$ by
$$\epsilon: ((v_1,k_1)...(v_r,k_r)) \mapsto e(v_1,k_1)...e(v_r,k_r).$$

\begin{prop}  Let $\Gamma$ be a generalized layered graph with a unique
minimal vertex. Then $\epsilon(\mathbf B(\Gamma))$ spans $A(\Gamma)$.
\end{prop}

\begin{proof}  Let $*$ denote the unique minimal vertex of $\Gamma$ and let $e \in E$.
If $h(e) \ne *$ then $1 + \sum_{j=1}^{l(e)} (-1)^j a_j(e)t^j = P_e(t) = P_{t(e)}(t)P_{h(e)}(t)^{-1}$
and so each $a_i(e)$ is in the subalgebra of $A(\Gamma)$ generated by the
coefficients of $P_{t(e)}(t)$ and $P_{h(e)}(t)$, and hence in the
subalgebra generated by all $e(v,j), v \in V_+, 1 \le j \le |v|$.
If $h(e) = *$ then $P_e(t) = P_{t(e)}(t)$ and so $a_i(e) = e(t(e),i)$
for $1 \le i \le l(e)$.  Thus $A(\Gamma)$ is generated by
$\{e(v,i)| v \in V_+, 1\le i \le |v|\}$.  Therefore $A(\Gamma)$
is spanned by the set of all products of these elements. Proposition 2.2 then gives the result.
\end{proof}

\section{Operations on graphs}

We define several operations on graphs (adding a vertex,
adding an edge, inversion, formation of bouquets) and
discuss the relations between the corresponding algebras.
As a consequence of our results on adding a vertex we see 
(Theorem 3.2) that the spanning set of Proposition 2.3 is a linear basis for $A(\Gamma)$.

Let $\Gamma = (G,|.|)$, where $G = (V,E)$, be a generalized layered graph,
$e \in E$ be an edge of length greater than $1$, and
$i$ an integer, $0 < i < l(e)$.  We define a new graph $G^w = (V^w,E^w)$ by
placing a new vertex, $w$, on the edge $e$.
Thus $$V^w = V \cup \{w\}$$
and
$$E^w = (E\setminus \{e\}) \cup \{e_1,e_2\}$$
where $$t(e_1) = t(e), h(e_1) = t(e_2)=w, h(e_2) = h(e).$$
Now extend the rank function $|.|$ to $ V^w$ by $|w| = |h(e_2)| + i$.  
Then $\Gamma^w = (G^w,|.|)$ is a generalized layered graph and $$ l(e_1) = l(e) - i, l(e_2) = i.$$
Define a map
$$\tilde\iota: E^{\sharp} \rightarrow T(E^{w{\sharp}})$$
by
$$\tilde\iota : a_j(f) \mapsto a_j(f)$$
for all $f\in E, f \ne e, 1 \le j \le l(f),$ and
$$\tilde\iota: a_j(e) \mapsto \sum_{\max (0,i+j-l(e_1)) \leq k \leq \min (i,j)}a_{j-k}(e_1)a_k(e_2)$$
for $1 \le j \le l(e).$
The map $\tilde\iota$ extends to  homomorphisms of graded algebras
$$\tilde\iota: T(E^{\sharp}) \rightarrow T(E^{w\sharp})$$ and
$$\tilde\iota: T(E^{\sharp})[t] \rightarrow T(E^{w\sharp})[t]$$
and we have
$$\tilde\iota(P_e(t)) = P_{e_1}(t)P_{e_2}(t).$$

Let $\pi = (f_1,...,f_r)$ be a path in $\Gamma$.
If $f_1,...,f_r \ne e$ define
$\tilde\iota(\pi)$ to be the path $(f_1,...,f_r)$ in $\Gamma^w$.
If $f_j = e$ define $\tilde\iota(\pi)$ to be the path
$(f_1,...,f_{j-1},e_1,e_2,f_{j+1},...,f_r)$ in $\Gamma^w$.  Then we have
$$\tilde\iota(P_{\pi}(t)) = P_{\tilde\iota(\pi)}(t)$$
for any path $\pi$ in $E$.  If $\pi_1, \pi_2$ are paths in $\Gamma$
with $t(e_1) = t(e_2), h(e_1) = h(e_2)$,  then
$\tilde\iota\pi_1, \tilde\iota\pi_2$ are paths in $\Gamma^w$
with $t(\tilde\iota e_1) = t(\tilde\iota e_2), h(\tilde\iota e_1) = h(\tilde\iota e_2)$.
Thus the generators of the ideal $R$ in $T(E^{\sharp})$ are mapped
by $\tilde\iota$ to generators of the ideal $R^w$ in $T(E^{w\sharp})$
and so $\tilde\iota$ induces a homomorphism
$$\iota: A(\Gamma) \rightarrow A(\Gamma^w).$$

For $v \in V^w$, $1\le j\le |v|$ define $e^w(v,j)$ by $$P_v(t) = \sum_{j=0}^{|v|} (-1)^je^w(v,j)t^j.$$

\begin{lem}  (a) If $v \in V_+$, then $\iota(e(v,j)) = e^w(v,j).$

(b) $\iota(\epsilon (\mathbf B(\Gamma)) \subseteq \epsilon (\mathbf B(\Gamma^w)).$
\end{lem}

\begin{proof}  Part (a) follows from the definition of $e^w(v,j)$.  Part (b) then follows since
$e(v_1,j_1)\dots e(v_r,j_r) \in \epsilon (\mathbf B(\Gamma))$ if and only if 
$e^w(v_1,j_1)...e^w(v_r,j_r) \in \epsilon (\mathbf B(\Gamma^w)).$
\end{proof}

\begin{thm}  Let $\Gamma$ have a unique minimal vertex.  Then

(a) $\epsilon (\mathbf B(\Gamma))$ is a basis for $A(\Gamma);$

(b)  $\iota: A(\Gamma) \rightarrow A(\Gamma^w)$ is an injection.
\end{thm}

\begin{proof}  Set $s(\Gamma) = \sum_{e \in E} (l(e) - 1).$
Then $s(\Gamma) = 0$ if and only if every edge of $\Gamma$
has length one, i.e., if and only if $\Gamma$ is a layered graph
(in the sense of \cite{GRSW}). Then, by Theorem 4.3 of \cite{GRSW},
$\epsilon (\mathbf B(\Gamma))$ is a basis for $A(\Gamma)$ if $s(\Gamma) = 0$.

We proceed by induction on $s(\Gamma)$, assuming that $\epsilon (\mathbf B(\Gamma'))$
is a basis for $A(\Gamma')$ whenever
$s(\Gamma') < s(\Gamma)$.  Clearly $s(\Gamma^w) = s(\Gamma) - 1$
and so we have that $\epsilon (\mathbf B(\Gamma^w))$ is linearly independent.
Then, by Lemma 3.1, $\iota$ maps $\epsilon (\mathbf B(\Gamma))$ into the linearly
independent set $\epsilon (\mathbf B(\Gamma^w))$ in $A(\Gamma^w)$ and hence 
$\epsilon (\mathbf B(\Gamma))$
is linearly independent.  In view of Proposition 2.3, $\epsilon (\mathbf B(\Gamma))$
is a basis for $A(\Gamma)$ and $\iota$ is injective.
\end{proof}

Let $\Gamma = (G,|.|)$ where $G = (V,E)$ be a generalized layered graph and
$v,w \in V$ be vertices with $|v| > |w|$.
Recall (from Section 1) the definition of the graph obtained by adding a new edge from $v$ to $w$:
$G^e = (V,E^e), E^e = E \cup\{e\}$, t(e) = v, h(e) = w, $\Gamma^e = (G^e,|.|).$
Then $i^e = (i_V,i^e_E)$, where $i_V$ is the identity map on $V$ 
and $i^e_E$ is the injection of $E$ into $E^e$, is a morphism 
of $\Gamma$ to $\Gamma^e$ and so there is a corresponding 
homomorphism $A(i^e): A(\Gamma) \rightarrow A(\Gamma^e).$

\begin{prop}  Let $\Gamma$ have a unique minimal vertex.  Then
$A(i^e)$
is surjective.
\end{prop}

\begin{proof}  Since $\Gamma$ has a unique minimal vertex, $h(\pi_v) = h(\pi_w)$. Therefore $P_v(t) = P_e(t)P_w(t)$ in $A(\Gamma^e)[t]$ and so
$$P_e(t)=1 + \sum_{j=1}^{l(e)} (-1)^ja_j(e)t^j = P_v(t)P_w(t)^{-1}.$$
Since the coefficients of $P_v(t)$ and $P_w(t)$ are in the image of $A(i^e)$ we have that
$E^{e\sharp}$ is contained in the image of $A(i^e)$, giving the result.
\end{proof}

We will later (Corollary 5.3) see that if $v > w$ in $G$ then $A(i^e) $ is an isomorphism.

Let $\Gamma = (V,E)$ be a generalized layered graph with $V =  \coprod_{i=0}^n V_i.$
We define the inverted graph $\check\Gamma = (\check V, \check E)$
by reversing all edges.  Thus
$$\check V = \coprod_{i=0}^n \check V^i$$
where $$\check V^i = V_{n-i}$$ and
$$\check E = \{\check e|e \in E\}$$
where
$$t(\check e) = h(e), h(\check e) = t(e).$$
Note that $l(\check e) = l(e).$
Define $$\tilde\eta: E^{\sharp} \rightarrow \check E^{\sharp}$$ by
$$\tilde\eta: a_i(e) \rightarrow a_i(\check e).$$
Then  $\tilde\eta$ extends to  anti-isomorphisms of graded algebras
$$\tilde\eta: T(E^{\sharp}) \rightarrow T(\check E^{\sharp})$$ and
$$\tilde\eta: T(E^{\sharp})[t] \rightarrow T(\check E^{\sharp})[t].$$

If $\pi = (e_1,...,e_r)$ is a path in $\Gamma$ we set
$\tilde\eta(\pi) = (\check{e_r},...,\check{e_1})$,
a path in $\check\Gamma$.  Then
$$\tilde\eta(P_{\pi}[t]) = P_{\tilde\eta(\pi)}[t]$$
and so
$$\tilde\eta(R) = \check R.$$
Thus $\tilde\eta$ induces a anti-isomorphism
$$\eta:A(\Gamma) \rightarrow A(\check\Gamma).$$

Finally, let $\Gamma_1$ and $\Gamma_2$ be generalized layered
graphs with unique minimal vertices $*_1$ and $*_2$.
We define $\Gamma_1 \vee \Gamma_2$ to be the ``bouquet" obtained by
identifying the minimal vertices.  Clearly, $A(\Gamma_1 \vee \Gamma_2) $
is the free product of $A(\Gamma_1)$ and $A(\Gamma_2).$

\section{Hilbert series}

We will derive an expression for $H(A(\Gamma),z)$, the Hilbert
series of $A(\Gamma)$, when $\Gamma$ is a generalized
layered graph with a unique minimal element.
While this expression is identical to that of \cite{RSW2} for
layered graphs, we will present it in a form which is better suited to applications.

For convenience we write $h(z) = H(A(\Gamma),z)$,
where $\Gamma$ is a generalized layered  graph
with unique minimal element $*$ of level $0$. 
For $\mathbf b =((v_1,k_1),...,(v_r,k_r)) \in \mathbf B_1(\Gamma)$ define
$$|\mathbf b| = \sum_{s=1}^r k_s$$ and set
$$\mathbf B_1(\Gamma)_{[j]} = \{\mathbf b \in \mathbf B_1(\Gamma)|\  |\mathbf b| = j \}.$$
For any subset $X \subseteq \mathbf B_1(\Gamma)$ set
$$X_{[j]} = \mathbf B_1(\Gamma)_{[j]} \cap X$$
and $$||X|| = \sum_{j \ge 0} |X_{[j]}|z^j.$$
Now $\epsilon(\mathbf B_1(\Gamma)_{[j}]) \subseteq A(\Gamma)_{[j]}$ and
$\epsilon(\mathbf B(\Gamma)_{[j]})$ is a basis for $A(\Gamma)_{[j]}.$
If $\mathbf b = ((v_1,k_1),...,(v_r,k_r)), \mathbf c = ((w_1,l_1),...,(w_s,l_s)) \in \mathbf B_1(\Gamma)$ define
$\mathbf b \circ \mathbf c = ((v_1,k_1,...,(v_r,k_r),(w_1,l_1),...,(w_s,l_s))$.
For $v \in V_+$, define
$$\mathbf C_v(\Gamma) = \cup_{k=1}^{|v|} (v,k) \circ \mathbf B(\Gamma),$$
$$\mathbf B_v(\Gamma) = \mathbf C_v(\Gamma) \cap \mathbf B(\Gamma),$$
and
$$\mathbf D_v(\Gamma) =\mathbf C_v(\Gamma) \setminus \mathbf B_v(\Gamma).$$
Then $\mathbf B(\Gamma) = \{\emptyset\} \cup \bigcup_{v \in V_+} \mathbf B_v(\Gamma).$  
Let $h_v(z) = ||\mathbf B_v(\Gamma)||$.  Then
$$h(z) = ||\mathbf B(\Gamma)|| = 1 + \sum_{v  \in V_+}h_v(z). $$
Now  $$||\mathbf C_v|| =(z + \dots + z^{|v|})h(z) = z\left (\frac{z^{|v|} - 1}{z-1}\right
)h(z).$$
Since
$$\mathbf D_v = \{((v,k),(v_1,k_1)\dots (v_r,k_r))\ |\ 1 \le k \le
|v|, (v,k) \gtrdot $$ $$(v_1,k_1), ((v_1,k_1)\dots (v_l,k_l)) \in
\mathbf B(\Gamma)\}  $$ we have
$$\mathbf D_v = \bigcup_{v > v_1 > *} (v,|v|-|v_1|)\circ \mathbf B_{v_1}(\Gamma).$$
Then $||\mathbf D_v|| = \sum_{v > v_1 > *} z^{|v|-|v_1|}h_{v_1}(z)$
and so
$$h_v(z) =
z\left (\frac{z^{|v|} - 1}{z-1}\right )h(z) - \sum_{v > w>*} z^{|v| -
|w|}h_w(z).$$

This equation may be written in matrix form.  Arrange the elements
of $V_+ $ in increasing order and index the elements
of vectors and matrices by this ordered set.  Let ${\mathbf h}(z)$
denote the column vector with entry $h_v(z)$ in the $v$-position, let
$\mathbf s$ denote the column vector with entry $\dfrac{z({z^{|v|} - 1 )}}{z-1}$
in the $v$-position,
let ${\mathbf 1}$ denote the column vector all of whose entries are $1$, and let $\zeta(z)$ denote
the matrix with entries $\zeta_{v,w}(z)$ for $ v,w \in V_+$ where
$\zeta_{v,w}(z) = z^{|v|-|w|}$ if $v \ge w$ and $0$ otherwise.
Then we have
$$\zeta(z){\mathbf h}(z) = {\mathbf s}h(z).$$

Now $N(z) = \zeta(z) - I$ is a strictly lower triangular matrix and so
$\zeta(z)$ is invertible.  In fact, $\zeta(z)^{-1} = I - N(z) + N(z)^2 -
\dots . $ For $v,w \in V_+$ set
$$\mu(v,w) = \sum_{v = v_1 > \dots > v_l = w}(-1)^{l+1},$$
the well-known M\"obius function of the partially ordered set $V_+$ (see Ch. 3 in \cite{S}).
Then the $(v,w)$-entry of $\zeta(z)^{-1}$ is
$$\mu(v,w)z^{|v|-|w|}.$$
Now
$$1 - h(z) = -{\mathbf 1}^T{\mathbf h}(z) = -{\mathbf 1}^T\zeta(z)^{-1}{\mathbf s}h(z).$$

Solving for $h(z)$ gives
$$h(z) = \frac{1}{1 - {\mathbf 1}^T\zeta(z)^{-1}{\mathbf s}}.$$

Using the expression given above for $\zeta(z)^{-1}$, we obtain the following result.
\medskip
\begin{thm}
Let $\Gamma$ be a generalized layered graph with unique minimal element $*$ of level $0$ 
and $h(z) $ denote the Hilbert series of $A(\Gamma)$.  Then

$$h(z) = \frac{1-z}{1-z + \sum_{v_1 > v_2 \dots > v_l>*} (-1)^l
(z^{|v_1| - |v_l|+1}-z^{|v_1|+1})}.$$
\end{thm}

This result can be restated using the important definition
$${\mathcal M}(\Gamma)(z) = \sum_{v > w \ge *} \mu(v,w)z^{|v|-|w|}.$$

\begin{cor}
Let $\Gamma$ be a generalized layered graph with unique minimal element $*$ of level $0$ 
and $h(z) $ denote the Hilbert series of $A(\Gamma)$.  Then
$$h(z) = \frac{1-z}{1-z{\mathcal M}(\Gamma)}.$$
\end{cor}

\begin{proof}  We show

$$\frac{1-z}{h(z)} = 1 - z + \sum_{v_1 > v_2\dots > v_l>*} (-1)^l
(z^{|v_1|-|v_l|+1} - z^{|v_1|+1})$$

$$  = 1  + \sum_{v_1 > v_2 \dots > v_l\ge *} (-1)^l  z^{|v_1|-|v_l|+1} = 1-z{\mathcal M}(\Gamma)(z).$$

The first equality is immediate from the theorem and the second follows by writing
$$\sum_{v_1 > v_2 \dots > v_l\ge *} (-1)^l  z^{|v_1|-|v_l|+1}$$
as
$$\sum_{v_1 > v_2 \dots > v_l> *} (-1)^l  z^{|v_1|-|v_l|+1} + \sum_{v_1 >
v_2 \dots > v_l = *} (-1)^l  z^{|v_1|-|v_l|+1}.$$
\end{proof}

As noted in Example 1.6, if $\Gamma = (G,|.|)$ where $G$ is a 
rooted tree, then $A(\Gamma)$ is the free algebra $T(E^{\sharp})$.  
By applying Corollary 4.2 in this situation we recover the well-known 
expression for the Hilbert series of a free algebra.
\begin{cor} 
Let $\Gamma = (G,|.|)$ be a generalized layered
graph where $G = (V,E)$ is a rooted tree. 
Let  $m_j, 1 \le j \le r,$ denote the number of edges of length $\ge j$.  Then
$$\frac{1}{H(\Gamma,z)}  = 1 - \sum_{j=1}^r m_jz^j.$$
\end{cor}
\begin{proof} Let
$E = \{e_1,...,e_n\}$ and recall that
$l(e)$ denotes the length of an edge $e$.   Now let $v,w \in V, v > w$ and let
$S(v,w) = \{u \in V|v \ge u \ge w\}.$ Then $$\mu(v,w) =
\sum_{\emptyset \ne T \subseteq S(v,w)} (-1)^{|T|}.$$  Thus
$\mu(v,w) = 1$ if $v = w, \mu(v,w) = -1$ if there is an edge $e$
from $v$ to $w$, and $\mu(v,w) = 0$ otherwise. Therefore
${\mathcal M}({\Gamma})(z) = |V| - \sum_{i=1}^n z^{l(e_i)}.$ Hence
$$1 - z{\mathcal M}(\Gamma)(z) = (1 - z) - z(|E| - \sum_{i=1}^nz^{l(e_i)}) = $$
$$(1-z) - z(\sum_{i=1}^n(1 - z^{l(e_i)}))
= (1-z)(1 -\sum_{i=1}^n (z + ... + z^{l(e_i)})) = $$ 
$$ (1-z)(1 - \sum_{j=1}^r m_jz^j).$$
\end{proof}

\section {Hilbert series and operations on graphs}

If $\Gamma$ has a unique minimal vertex $v_{\min}$, we define
$${\mathcal M}_{\circ}(\Gamma)(z) = \sum_{v \in V} \mu(v,v_{\min})z^{|v|-|v_{\min}|}.$$
Similarly, if $\Gamma$ has a unique maximal vertex $v_{\max}$ we
define
$${\mathcal M}^{\circ}(\Gamma)(z) = \sum_{v \in V} \mu(v_{\max},v)z^{|v_{\max}|-|v|}.$$
Let $e \in E$ be an edge of length greater than $1$ with $t(e) =
u$ and $ h(e) = v.$ Recall that for any $i, 1 \le i < l(e)$, we
have defined a new graph $\Gamma^w = (G^w,|.|)$ where $G^w = (V^w, E^w)$ by adding a new
vertex $w$ on the edge $e$ with $|w| - |v| = i.$ We define two
related graphs, $G^w_- = (V^w_-,E^w_-)$ and  $G^w_+ =
(V^w_+,E^w_+)$ by $V^w_- = \{v \in V^w|w \ge v\}, E^w_- = \{e \in
E^w |\ t(e), h(e) \in V^w_-\}, V^w_+ = \{v \in V^w|v \ge w\}, E^w_+
= \{e \in E^w |\ t(e), h(e) \in V^w_+\}.$ 
Then the restrictions of the rank function $|.|$ to $V^w_-$ and to 
$V^w_+$ (which we continue to denote by $|.|$) are rank functions.  
Thus we have generalized layered graphs $\Gamma^w_- = (G^w_-,|.|)$
and $\Gamma^w_+ = (G^w_+,|.|)$.
\begin{prop}  $${\mathcal M}(\Gamma^w)(z) = {\mathcal M}(\Gamma)(z) +
{\mathcal M}_{\circ}(\Gamma^w_+)(z) \cdot {\mathcal M}^{\circ}(\Gamma^w_-)(z).$$
\end{prop}

\begin{proof}  For $u', v' \in V^w , u' \ge w \ge v'$ define
$$\mu^w(u',v') = \mathop{\sum_{u' = x_1 > \dots > x_l = v'}}_{w \in \{x_1,...,x_l\}} (-1)^{l+1}.$$

Then $${\mathcal M}(\Gamma^w)(z) - {\mathcal M}(\Gamma)(z) =
\sum_{u' > w > v'} \mu^w(u',v')z^{|u'|-|v'|} $$
$$+ \sum_{u' > w} \mu^w(u',w)z^{|u'|-|w|}
+ \sum_{w > v'} \mu^w(w,v')z^{|w|-|v'|} + \mu^w(w,w).$$

Now 
$${\mathcal M}_{\circ}(\Gamma^w_+)(z) -1 = \sum_{u'>w}\mu^w(u',w)z^{|u'|-|w|} = -\sum_{u' \ge u'' \ge u} \mu(u',u'')z^{|u'|-|w|}  ,$$
and

$${\mathcal M}^{\circ}(\Gamma^w_-)(z) -1 = \sum_{w > v'} \mu^w(w,v')z^{|w|-|v'|} = -\sum_{v \ge v'' \ge v'} \mu(v'',v')z^{|w|-|v'|}.$$
Also, if $u' > w > v'$, we have

$$\mu^w(u',v')z^{|u'|-|v'|} = \sum_{u' \ge u'' \ge u > v  \ge v'' \ge v'} \mu(u',u'')\mu(v'',v')z^{|u'|-|v'|} = $$
$$(\sum_{u' \ge u'' \ge u}\mu(u',u'')z^{|u'|-|w|})(\sum_{v \ge v'' \ge v'}\mu(v'',v')z^{|w|-|v'|})= $$
$$ ({\mathcal M}_{\circ}(\Gamma^w_+)(z)-1)({\mathcal M}^{\circ}(\Gamma^w_-)(z)-1),$$
giving the result.
\end{proof}

Now let $\Gamma = (V,E), v,w \in V, |v| > |w|$.  Recall that we have defined a graph
$\Gamma^e = (V^e,E^e)$ by adjoining an edge $e$ to $E$ with $t(e) = v, h(e) = w$. 
If $a,b \in V$ write $a > b $ if there is a path in $E$ from $a$ to $b$.Then the following proposition is immediate from the
definition of $\mathcal M$.
\medskip
\begin{prop}$${\mathcal M}(\Gamma^e)(z) - {\mathcal M}(\Gamma)(z) =
\sum (-1)^{l+m} z^{|v_1|-|w_m|}$$
where the sum is taken over all sequences $v_1 > ... > v_l \ge v,$ $ w \ge w_1 > ... > w_m,$  
such that  there is no path in $E$ from $v_l$ to $w_1$.
\end{prop}

\begin{cor}  If $\Gamma^e$ is obtained from $\Gamma$ by adjoining an edge from $v$ to $w$ where there is a path in $\Gamma$ from $v$ to $w$  then
${\mathcal M}(\Gamma^e)(z) = {\mathcal M}(\Gamma)(z)$. Hence
$H(A(\Gamma^e),z) = H(A(\Gamma),z)$ so $A(i^e):A(\Gamma) \rightarrow A(\Gamma^e)$ is an isomorphism.
\end{cor}

\begin{proof}  If there is a path from $v$ to $w$ in $\Gamma$ then the sum occuring in the proposition is vacuous.
\end{proof}

Next let $\Gamma = (G,|.|)$ where $G = (V, E)$ have unique maximal vertex and a unique minimal vertex.  We have defined the inverted graph $\check\Gamma$.
Since $v_1 > ... > v_l$ in $V$ if and only if $v_l > ... > v_1$ in $\check V$, the following proposition is immediate from the
definition of $\mathcal M$.

\begin{prop} Let $\Gamma$ have a unique maximal vertex and a unique minimal vertex.
Then ${\mathcal M}(\check\Gamma)(z) = {\mathcal M}(\Gamma)(z)$ and so $H(A(\check\Gamma),z) = H(A(\Gamma),z).$
\end{prop}

If  $\Gamma_1$ and $\Gamma_2$ are generalized layered
graphs with unique minimal vertices $*_1$ and $*_2$, we have
defined $\Gamma_1 \vee \Gamma_2$ to be the ``bouquet" obtained by
identifying the minimal vertices. The following proposition is clear.

\begin{prop}  ${\mathcal M}(\Gamma_1 \vee \Gamma_2)(z) = {\mathcal M}(\Gamma_1)(z) + {\mathcal M}(\Gamma_2)(z) - 1.$
\end{prop}

We remark that this implies
$$\frac{1}{H(A(\Gamma_1 \vee \Gamma_2),z)} = \frac{1}{H(A(\Gamma_1),z)} + \frac{1}{H(A( \Gamma_2),z)} - 1.$$
Of course, this already followed from our previous observation that $A(\Gamma_1 \vee  \Gamma_2)$
is the free product of $A(\Gamma_1)$ and $A(\Gamma_2).$

Finally, let $\Gamma_1$ and $\Gamma_2$ be generalized layered
graphs with unique minimal vertices $v_{\min,1}$ and $v_{\min,2}$
of level $0$ and unique maximal vertices $v_{\max,1}$ and
$v_{\max,2}$ of level $d$.  We define $\Gamma_1 \diamondsuit
\Gamma_2$ to be the ``double bouquet"  obtained by identifying
$v_{\min,1} $ with $v_{\min,2}$ and identifying $v_{\max,1}$ with
$v_{\max,2}.$

\begin{prop}  ${\mathcal M}(\Gamma_1 \diamondsuit \Gamma_2)(z) = {\mathcal M}(\Gamma_1)(z) + {\mathcal M}(\Gamma_2)(z) +
2 - z^d$.
\end{prop}

\begin{proof}  Let $S_i$ denote the set of sequences $u_1 > ... > u_l$ in $V_i$ for $i = 1,2$
and let $S_{12}$ denote the set of sequences $u_1 > ... > u_l$ in
$V_1 \diamondsuit V_2$.  Write $v_{\max} = v_{\max,1} =
v_{\max,2}$ and $v_{\min} = v_{\min,1} = v_{\min,2}$ in $V_1
\diamondsuit V_2$. Clearly $S_{12} = S_1 \cup S_2$ and $S_1 \cap
S_2 = \{(v_{\max}), (v_{\min}),(v_{\max},v_{\min})\}$, giving the
result.
\end{proof}

Let $\Delta(d)$ denote the generalized layered graph with two
vertices $max$ and $min$ of levels $d$ and $0$ respectively and
with one edge $e$ with $t(e) = max, h(e) = min.$

\begin{cor} $$\frac{1}{H(A(\Gamma_1 \diamondsuit \Gamma_2),z)} = \frac{1}{H(A(\Gamma_1),z)} +
\frac{1}{H(A(\Gamma_2),z)} - \frac{1}{H(A(\Delta(d),z)}.$$
\end{cor}

Recall that according to \cite{Golod} if $V$ is a graded vector
space with graded dimension $H(V,z)$ and $R$ is a graded subspace
of the tensor algebra $T(V)$ with graded dimension $H(R,z)$, the
quotient $A = T(V)/<R>$ is said to be a noncommutative complete
intersection if $$H(A,z) = \frac {1}{1 - H(V,z) + H(R,z)}.$$

\begin{cor} Let $\Gamma_1$ and $\Gamma_2$ be generalized layered graphs
with unique minimal vertices $v_{\min,1}$ and $v_{\min,2}$ of
level $0$ and unique maximal vertices $v_{\max,1}$ and
$v_{\max,2}$ of level $d$.  Assume that $A(\Gamma_1)$ and
$A(\Gamma_2)$ are noncommutative complete intersections.  Then
$A(\Gamma_1 \diamondsuit \Gamma_2)$ is a noncommutative complete
intersection.
\end{cor}

\begin{proof}  Let $G_i$ denote space of generators for $A(\Gamma_i)$
and $R_i$ denote the space of relations.  Since $A(\Delta(d))$ is
a free algebra on generators of degrees $1,2,...,d$ we have
$H(A(\Delta(d),z)^{-1} = 1 - z -  z^2 - ... - z^d$.  Thus
$$\frac{1}{H(A(\Gamma_1 \diamondsuit \Gamma_2),z)} = $$
$$1 - H(G_1,z) + H(R_1,z) + 1 - H(G_2,z) + H(R_2,z) -1 + z + ... + z^d = $$
$$1 - H(G_1 + G_2,z) + H(R_1 + R_2,z) + z + ... + z^d.$$
Now the generators for $\Gamma_1 \diamondsuit \Gamma_2 $
are just the generators for $\Gamma_1$ and for $\Gamma_2$ and
the relations for $\Gamma_1 \diamondsuit \Gamma_2$ are just the
relations for $\Gamma_1$ and for $\Gamma_2$ together with the relations
stating that $P_{\pi_1}(t) = P_{\pi_2}(t)$ where $\pi_i$ is a path from
$v_{max,i}$ to $v_{min,i}$.  Since $P_{\pi_1}(t)$ has degree $d$,
there is one such relation of degree $j$ for $j = 1,...,d.$
Thus the space of relations for $\Gamma_1 \diamondsuit \Gamma_2$ has graded
dimension $H(R_1 + R_2,z) + z + ... + z^d$ and the corollary is proved.
\end{proof}

\begin{example}. Using Corollaries 4.3 and 5.7 we see that
the Hilbert series for the algebra $A(\Gamma([3,(12)])$ (defined
in Example 1.7) is $\frac{1}{1 - 3z -z^2 + z^3}$.
\end{example}

\enddocument